\documentclass{article}
\usepackage[utf8]{inputenc}
\usepackage{hyperref}
\usepackage{amsmath, amsthm, amssymb, latexsym,epsfig,amsthm,enumerate,multicol,wasysym}
\usepackage{xcolor, comment}
\usepackage{bbm}
\usepackage{mathtools}
\usepackage{enumitem}
\usepackage[pagewise]{lineno}
\usepackage[symbol]{footmisc}

\newcommand{\regret}{\textsc{Regret}}
\newcommand{\R}{\mathbb{R}}
\newcommand{\cost}{\textsc{Cost}}
\newcommand{\ecost}{\textsc{ECost}}

\newcommand{\opt}{\text{opt}}
\newcommand{\prior}{\text{Prior}}
\newcommand{\eprior}{\emph{Prior}}
\newcommand{\mx}{\text{max}}
\newcommand{\emx}{\emph{max}}
\newcommand{\tame}{\text{TAME}}
\newcommand{\etame}{\emph{TAME}}
\newcommand{\bayes}{\text{Bayes}}
\newcommand{\ebayes}{\emph{Bayes}}
\newcommand{\E}{\text{E}}
\newcommand{\ag}{\text{Ag}}
\newcommand{\cK}{\mathcal{K}}
\newcommand\extrafootertext[1]{%
    \bgroup
    \renewcommand\thefootnote{\fnsymbol{footnote}}%
    \renewcommand\thempfootnote{\fnsymbol{mpfootnote}}%
    \footnotetext[0]{#1}%
    \egroup
}

\newtheorem*{pde}{PDE Assumption}
\newtheorem{thm}{Theorem}
\newtheorem{lem}{Lemma}

\title{Almost Optimal Agnostic Control of Unknown Linear Dynamics}
\author{J. Carruth, M. Eggl, C. Fefferman, C. Rowley}
\date{January 2024}

\begin{document}

\maketitle

\begin{abstract}
    We consider a simple control problem in which the underlying dynamics depend on a parameter $a$ that is unknown and must be learned. We study three variants of the control problem: Bayesian control, in which we have a prior belief about $a$; bounded agnostic control, in which we have no prior belief about $a$ but we assume that $a$ belongs to a bounded set; and fully agnostic control, in which $a$ is allowed to be an arbitrary real number about which we have no prior belief. In the Bayesian variant, a control strategy is optimal if it minimizes a certain expected cost. In the agnostic variants, a control strategy is optimal if it minimizes a quantity called the worst-case regret. For the Bayesian and bounded agnostic variants above, we produce optimal control strategies. For the fully agnostic variant, we produce almost optimal control strategies, i.e., for any $\varepsilon>0$ we produce a strategy that minimizes the worst-case regret to within a multiplicative factor of $(1+\varepsilon)$.
\end{abstract}

The purpose of this note is to announce the results of our companion papers \cite{almostoptimal2023, bounded2023}. These papers explore a new flavor of adaptive control theory, which we call ``agnostic control''; see also \cite{carruth2023bounded, carruth2022controlling, Fefferman:2021, gurevich2022optimal}. While our exposition here borrows heavily from the introductions of \cite{almostoptimal2023,bounded2023}, we think the results benefit from a unified presentation. Moreover, we give here a more detailed overview of the results of \cite{bounded2023} than is given in the introduction to that paper.\extrafootertext{This work was supported by AFOSR grant FA9550-19-1-0005 and by the Joachim Herz Foundation..}

Many works in adaptive control theory attempt to control a system whose underlying dynamics are initially unknown and must be learned from observation. The goal is then to bound $\regret$, a quantity defined by comparing our expected cost with that incurred by an opponent who knows the underlying dynamics and plays optimally. Typically one tries to achieve a regret whose order of magnitude is as small as possible after a long time. Adaptive control theory has extensive practical applications; see, e.g., \cite{bertsekas2012dynamic, cesa2006prediction, hazancontrol, powell2007approximate} for some examples.

In some applications, we don't have the luxury of waiting for a long time. This is the case, e.g., for a pilot attempting to land an airplane following the sudden loss of a wing, as in \cite{Brazy:2009}. Our goal here is to achieve the absolute minimum possible regret over a fixed, finite time horizon. This objective poses formidable mathematical challenges, even for simple model systems.

We will study a one-dimensional, linear model system whose dynamics depend on a single unknown parameter $a$. When $a$ is large positive, the system is highly unstable. (There is no ``stabilizing gain'' for all $a$.) We will make progressively weaker assumptions about the unknown parameter $a$---eventually, we will assume that $a$ may be any real number and we won't assume that we are given a Bayesian prior probability distribution for it.

We now give a precise statement of our problem.

\section*{The Model System}

Our system consists of a particle moving in one dimension, influenced by our control and buffeted by noise. The position of our particle at time $t$ is denoted by $q(t) \in \R$. At each time $t$, we may specify a ``control'' $u(t) \in \R$, determined by history up to time $t$, i.e., by $(q(s))_{s \in [0,t]}$. A ``strategy'' (aka ``policy'') is a rule for specifying $u(t)$ in terms of $(q(s))_{s \in [0,t]}$ for each $t$. We write $\sigma, \sigma', \sigma^*, \text{etc.}$ to denote strategies. The noise is provided by a standard Brownian motion $(W(t))_{t\ge 0}$.

The particle moves according to the stochastic ODE
\begin{equation}\label{eq: intro 1}
dq(t) = \big(aq(t) + u(t)\big)dt + dW(t), \qquad q(0) = q_0,
\end{equation}
where $a$ and $q_0$ are real parameters. Due to the noise in \eqref{eq: intro 1}, $q(t)$ and $u(t)$ are random variables; these random variables depend on our strategy $\sigma$, and we often write $q^\sigma(t)$, $u^\sigma(t)$ to make that dependence explicit.

Over a time horizon $T>0$, we incur a $\cost$, given\footnote[2]{By rescaling, we can consider seemingly different cost functions of the form $\int_0^T(q^2+\lambda u^2)$ for $\lambda >0$.} by
\begin{equation}\label{eq: intro 2}
    \cost(\sigma,a) = \int_0^T \big\{ (q^\sigma(t))^2 + (u^\sigma(t))^2\big\} dt.
\end{equation}
This quantity is a random variable determined by $a, q_0,T$ and our strategy $\sigma$. Here, the starting position $q_0$ and time horizon $T$ are fixed and known. 

We would like to pick our strategy $\sigma$ to keep our cost as low as possible. We examine several variants of the above control problem, making successively weaker assumptions regarding our knowledge of the parameter $a$. The first variant is simply the classical case, in which $a$ is a known real number. In the second variant, we assume that the parameter $a$ is unknown, but subject to a given prior probability distribution supported on a bounded interval. In the third variant, we assume that the parameter $a$ belongs to a bounded interval, but is otherwise unknown (in particular, we do \emph{not} assume that we are given a prior belief about $a$). In the fourth and final variant, we assume that $a$ is unknown and may be any real number (again, we do not assume that we are given a prior belief about $a$). We refer to the third and fourth variants, in which we are not given a prior belief about $a$, as \emph{agnostic control}.

\subsection*{Variant I: Classical Control}
We suppose first that the parameter $a$ is known. We write $\ecost(\sigma, a; T,q_0)$, or sometimes $\ecost(\sigma,a)$, to denote the expected $\cost$ incurred by executing a given strategy $\sigma$. Our task is to pick $\sigma$ to minimize $\ecost(\sigma, a; T,q_0)$. As shown in textbooks (e.g., \cite{astrom}), there is an elementary formula for the optimal strategy, denoted $\sigma_\opt(a)$, given by
\[
u(t) =- \kappa(T-t,a)q(t),
\]
where
\[
\kappa(s,a) = \frac{\tanh(s\sqrt{a^2+1})}{\sqrt{a^2+1}-a\tanh(s\sqrt{a^2+1})}.
\]
We refer to $\sigma_\opt(a)$ as the \emph{optimal known}-$a$ \emph{strategy}. It will be important later to note that $\sigma_\opt(a)$ satisfies the inequality
\begin{equation}\label{eq: 3}
    |u(t)|\le C \max\{a,1\}\cdot |q(t)|\;\text{for an absolute constant } C.
\end{equation}

\subsection*{Variant II: Bayesian Control}
We now suppose that the parameter $a$ is unknown, but is subject to a given prior probability distribution $d\prior(a)$ supported in an interval $[-a_\mx, a_\mx]$. Our goal is then to pick a strategy $\sigma$ to minimize our expected cost, given by
\begin{equation}\label{eq: 4}
\ecost(\sigma, d\prior) = \int_{-a_\mx}^{a_\mx} \ecost(\sigma, a) \ d\prior(a).
\end{equation}

Before presenting rigorous results, we provide a heuristic discussion. 

First of all, since $d\prior$ is supported in $[-a_\mx, a_\mx]$, a glance at \eqref{eq: 3} suggests that our optimal strategy $\sigma$ will satisfy
\begin{equation}\label{eq: 4.5}
|u^\sigma(t)| \le C a_\mx |q^\sigma(t)|.
\end{equation}
 In \cite{bounded2023}, we introduce the notion of a \emph{tame strategy} $\sigma$, which satisfies the estimate
\begin{equation}\label{eq: 5}
|u^\sigma(t)| \le C_{\tame}^\sigma [|q^\sigma(t)| + 1] \quad (\text{for all}\; t \in [0,T])
\end{equation}
with probability 1, for a constant $C_{\tame}^\sigma$ called a \emph{tame constant} for $\sigma$ (note that $C_\tame$ may depend on $a_\mx$). Thus, we expect that the optimal strategy for Bayesian control will be tame.

Next, we note a major simplification. In principle, a strategy $\sigma$ is a one-parameter family of functions on an infinite-dimensional space, because for each $t$ it specifies $u(t)$ in terms of the path $(q(s))_{s \in [0,t]}$. However, reasoning heuristically, one computes that the posterior probability distribution for the unknown $a$, given a past history $(q(s))_{s \in [0,t]}$ is determined by the prior $d\prior(a)$, together with the two observable quantities
\begin{equation}
\zeta_1(t) = \int_0^t q(s)[dq(s) - u(s) ds]\quad \text{and}\quad \zeta_2(t) = \int_0^t (q(s))^2\ ds \ge 0.
\end{equation}
Therefore, it is natural to suppose that the optimal strategy $\sigma_{\bayes}(d\prior)$ takes the form
\begin{equation}\label{eq: 7}
    u(t) = \tilde{u}(q(t), t, \zeta_1(t), \zeta_2(t))
\end{equation}
for a function $\tilde{u}$ on $\R \times [0,T] \times \R \times [0,\infty)$.

So, instead of looking for a one-parameter family of functions on an infinite-dimensional space, we merely have to specify a function $\tilde{u}$ of four variables.

It isn't hard to apply heuristic reasoning to derive a PDE for the function $\tilde{u}$ in \eqref{eq: 7}. To do so, we introduce the \emph{cost-to-go}, $S(q,t,\zeta_1,\zeta_2)$, defined as the expected value of
\begin{equation}\label{eq: 8}
    \int_t^T \{(q^\sigma(t))^2 + (u^\sigma(t))^2\} \ ds,
\end{equation}
conditioned on
\begin{equation}\label{eq: 9}
    q(t) = q, \quad \zeta_1(t) = \zeta_1, \quad \zeta_2(t) = \zeta_2,
\end{equation}
with our strategy $\sigma$ picked to minimize \eqref{eq: 8}. Clearly,
\begin{equation}\label{eq: 10}
\begin{split}
S(&q,t,\zeta_1,\zeta_2) \\
&= \min_u\{(q^2 + u^2)dt + \E[S(q+dq, t+dt, \zeta_1 + d\zeta_1, \zeta_2 + d\zeta_2)]\} + o(dt),
\end{split}
\end{equation}
where $(q,t,\zeta_1,\zeta_2)$ evolves to $(q+dq, t + dt, \zeta_1 + d\zeta_1, \zeta_2 + d\zeta_2)$ if we apply the control $u$ at time $t$. Here, $\E[\cdots]$ denotes expected value conditioned on \eqref{eq: 9}.

Moreover, the optimal control at time $t$ (given \eqref{eq: 9}) is precisely the value of $u$ that minimizes the right-hand side of \eqref{eq: 10}; we denote it $\tilde{u}(t)$.

Taylor-expanding the right-hand side of \eqref{eq: 10}, and taking $dt\rightarrow 0$, we arrive at the \emph{Bellman equation}
\begin{equation}\label{eq: 11}
    \begin{split}
        0 = &\partial_t S + (\bar{a}(\zeta_1,\zeta_2) q + \tilde{u})\partial_q S + \bar{a}(\zeta_1,\zeta_2) q^2 \partial_{\zeta_1}S + q^2 \partial_{\zeta_2}S+ \frac{1}{2}\partial_q^2 S \\& + q \partial_{q\zeta_1} S  + \frac{1}{2}q^2 \partial_{\zeta_1}^2 S + (q^2 + \tilde{u}^2),
    \end{split}
\end{equation}
where $\bar{a}(\zeta_1, \zeta_2)$ is the posterior expected value of $a$ given \eqref{eq: 9}; explicitly,
\begin{equation}\label{eq: 12}
    \bar{a}(\zeta_1,\zeta_2) = \frac{ \int_{-a_\mx}^{a_\mx}a \exp\big( - \frac{a^2}{2}\zeta_2 + a \zeta_1\big)\ d\prior(a) }{\int_{-a_\mx}^{a_\mx} \exp\big( - \frac{a^2}{2}\zeta_2 + a \zeta_1\big)\ d\prior(a) } .
\end{equation}
Moreover, the minimizer $\tilde{u}$ for the right-hand side of \eqref{eq: 10} is given by
\begin{equation}\label{eq: 13}
    \tilde{u}(q,t,\zeta_1,\zeta_2) = - \frac{1}{2} \partial_qS(q,t,\zeta_1,\zeta_2).
\end{equation}
Together with \eqref{eq: 11}, we impose the obvious terminal condition
\begin{equation}\label{eq: 14}
    S|_{t=T} = 0,
\end{equation}
and the natural requirement
\begin{equation}\label{eq: 15}
    S \ge 0.
\end{equation}

Our plan to solve for the optimal Bayesian control is thus to solve \eqref{eq: 11}--\eqref{eq: 15} for $S$ and $\tilde{u}$, and then set $u^{\sigma_{\bayes}(d\prior)}(t):= \tilde{u}(q(t), t, \zeta_1(t), \zeta_2(t))$.

We have produced numerical solutions to \eqref{eq: 11}--\eqref{eq: 15}, but we don't have rigorous proofs of existence or regularity. We proceed by imposing the following assumption.

\begin{pde}
    Equations \eqref{eq: 11}--\eqref{eq: 15} admit a solution $S \in C^{2,1}(\R\times[0,T]\times \R \times [0, \infty))$, satisfying the estimates
\begin{equation}\label{eq: 16}
    |\partial_{q,t,\zeta_1,\zeta_2}^\alpha S| \le K \cdot [1 + |q| + |\zeta_1| + \zeta_2|]^{m_0}\;\text{a.e. for}\; |\alpha| \le 3,
\end{equation}
and
\begin{equation}\label{eq: 17}
    |\tilde{u}| \le C_{\etame} \cdot [1 + |q|]\;\text{for all}\; (q, t, \zeta_1, \zeta_2),
\end{equation}
for some $K$, $m_0$, $C_\etame$.
\end{pde}

Assumption \eqref{eq: 17} asserts that our strategy $\sigma_{\bayes}(d\prior)$, given by \eqref{eq: 11}--\eqref{eq: 15}, is a tame strategy, as expected.

Our numerical simulations appear to confirm \eqref{eq: 16}, \eqref{eq: 17}. Accordingly, our \emph{PDE Assumption} seems safe.

We are ready to present our rigorous results on optimal Bayesian control; these are proved in \cite{bounded2023}.

\begin{thm}[Optimal Bayesian Strategy]\label{thm: 1}
Fix a probability distribution $d\eprior$, supported on $[-a_\emx, a_\emx]$, and suppose our \emph{PDE Assumption} is satisfied. Let $\sigma = \sigma_{\ebayes}(d\eprior)$ be the strategy obtained by solving \eqref{eq: 11}--\eqref{eq: 15}. Then
\begin{enumerate}
    \item[\emph{(A)}] $\ecost(\sigma, d\eprior) = S(q_0,0,0,0)$, with $S$ as in \eqref{eq: 11}--\eqref{eq: 15}.
    \item[\emph{(B)}] Let $\sigma'$ be any other strategy. Then
    \[
    \ecost(\sigma', d\eprior) \ge \ecost(\sigma, d\eprior),
    \]
    with equality only when we have
    \[
    u^{\sigma'}(t) = u^\sigma(t)\;\text{for a.e.}\; t\quad \text{and}\quad q^{\sigma'}(t) = q^\sigma(t)\;\text{for all}\; t,
    \]
    with probability 1.
\end{enumerate}
\end{thm}
When the competing strategy $\sigma'$ is assumed to be tame, we can sharpen the above uniqueness assertion (B) to a quantitative result.

\begin{thm}[Quantitative Uniqueness of the Optimal Bayesian Strategy]\label{thm: 2}
Let $d\eprior$ and $\sigma = \sigma_{\ebayes}(d\eprior)$ be as in Theorem \ref{thm: 1}. Given $\varepsilon>0$, and given a constant $\hat{C}$, there exists $\delta>0$ for which the following holds.

Let $\sigma'$ be a tame strategy with tame constant at most $\hat{C}$. If
\[
\ecost(\sigma', d\eprior) \le \ecost(\sigma,d\eprior) + \delta,
\]
then the expected value of
\[
\int_0^T \{ |q^\sigma(t) - q^{\sigma'}(t)|^2 + |u^\sigma(t) - u^{\sigma'}(t)|^2\} \ dt
\]
is less than $\varepsilon$.
\end{thm}

Theorem \ref{thm: 2} plays a crucial r\^{o}le in our analysis of agnostic control for bounded $a$ (see \cite{bounded2023} for details).

We now discuss an issue arising in the proofs of our results on Bayesian control: We need a rigorous definition of a strategy. Clearly, the phrase ``a rule for specifying $u(t)$ in terms of past history'' isn't precise.

We want to allow $u(t)$ to depend discontinuously on past history $(q(s))_{s \in [0,t]}$. For instance, we should be allowed to set
\[
u(t) = \begin{cases}
    -q(t) &\text{if}\; |q(t)| > 1,\\
    0 &\text{otherwise}.
\end{cases}
\]
On the other hand, we had better make sure that we can produce solutions of our stochastic ODE
\begin{equation}\label{eq: p0}
dq = (aq+u)dt + dW.
\end{equation}
Without the noise $dW$, we have a standard ODE, and the usual existence and uniqueness theorems for ODE would require Lipschitz continuity of $u$.

We proceed as follows.

At first we fix a partition
\begin{equation}\label{eq: p1}
    0 = t_0 < t_1 < \dots < t_N = T
\end{equation}
of the time interval $[0,T]$. We restrict ourselves to strategies $\sigma$ in which the control $u(t)$ is constant in each interval $[t_\nu, t_{\nu+1})$, and in which, for each $\nu$, $u(t_\nu)$ is determined by $(q(t_\gamma))_{\gamma \le \nu}$, together with ``coin flips'' $\vec{\xi} = (\xi_1, \xi_2, \dots) \in \{0,1\}^\mathbb{N}$. We assume that $u(t_\nu)$ is a Borel measurable function of $(q(t_1), \dots, q(t_\nu), \vec{\xi})$, and that for all $\nu$ we have
\[
|u(t_\nu)| \le C_{\tame} [ |q(t_\nu)| + 1].
\]
We call such a strategy a \emph{tame strategy associated to the partiton \eqref{eq: p1}} with \emph{tame constant} $C_\tame$. For such strategies, it is easy to define the solutions $q^\sigma(t)$, $u^\sigma(t)$ of our stochastic ODE \eqref{eq: p0}.

Most of our work lies in controlling and optimizing tame strategies associated to a sufficiently fine partition \eqref{eq: p1}. In particular, we prove approximate versions of Theorems \ref{thm: 1} and \ref{thm: 2} in the setting of such tame strategies.

We then define a tame strategy (not associated to any partition) by considering a sequence $\pi_1, \pi_2,\dots$ of ever-finer partitions of $[0,T]$. To each partition $\pi_n$ we associate a tame strategy $\sigma_n$ with a tame constant $C_\tame$ independent of $n$. If the resulting $q^{\sigma_n}(t)$ and $u^{\sigma_n}(t)$ tend to limits, in an appropriate sense, as $n \rightarrow \infty$, then we declare those limits $q(t)$, $u(t)$ to arise from a \emph{tame strategy}~$\sigma$.

Finally, we drop the restriction to tame strategies and consider general strategies. To do so, we consider a sequence $(\sigma_n)_{n = 1,2,\dots}$ of tame strategies, \emph{not} assumed to have a tame constant independent of $n$. If the relevant $q^{\sigma_n}(t)$ and $u^{\sigma_n}(t)$ converge, in a suitable sense, as $n \rightarrow \infty$, then we say that the limits $q(t)$ and $u(t)$ arise from a strategy $\sigma$.

It isn't hard to pass from tame strategies associated to partitions of $[0,T]$ to general tame strategies, and then to pass from such tame strategies to general strategies. The work in proving Theorems \ref{thm: 1} and \ref{thm: 2} lies in our close study of tame strategies associated to fine partitions. We refer the reader to \cite{bounded2023} for details.

\subsection*{Variant III: Agnostic Control for Bounded $a$}

We now suppose that our parameter $a$ is confined to a bounded interval $[-a_\mx, a_\mx]$ but is otherwise unknown. In particular, we don't assume that we are given a Bayesian prior probability distribution $d\prior(a)$. Consequently, we cannot define a notion of expected cost by formula \eqref{eq: 4}.

Instead, our goal will be to minimize \emph{worst-case regret}, defined by comparing the performance of our strategy with that of the optimal known-$a$ strategy $\sigma_\opt(a)$. We will introduce several variants of the notion of regret.

Let us fix a starting position $q_0$, a time horizon $T$, and an interval $[-a_\mx, a_\mx]$ guaranteed to contain the unknown $a$. To a given strategy $\sigma$, we associate the following functions on $[-a_\mx, a_\mx]$:
\begin{itemize}
    \item \underline{Additive Regret}, defined as 
    \[
    \text{AR}(\sigma,a) = \ecost(\sigma, a) - \ecost(\sigma_\opt(a),a)\ge 0.
    \]
    \item \underline{Multiplicative Regret} (aka ``competitive ratio''), defined as
    \[
    \text{MR}(\sigma, a) = \frac{\ecost(\sigma,a)}{\ecost(\sigma_\opt(a),a)}\ge 1.
    \]
    \item \underline{Hybrid Regret}, defined in terms of a parameter $\gamma >0$ by setting
    \[
    \text{HR}_\gamma(\sigma,a) = \frac{\ecost(\sigma,a)}{\ecost(\sigma_\opt(a),a)+\gamma}.
    \]
\end{itemize}
\sloppypar Writing $\regret(\sigma,a)$ to denote any one of the above three functions on $[-a_\mx, a_\mx]$, we define the \emph{worst-case regret}
\[
\regret^*(\sigma) = \sup  \big\{ \regret(\sigma, a) : a \in [-a_\mx, a_\mx ]\big\}.
\]
We seek a strategy $\sigma$ that minimizes worst-case regret.

The above notions are useful in different regimes. If we expect to pay a large cost, then we care more about multiplicative regret then about additive regret. (If we have to pay $10^9$ dollars, we are unimpressed by a savings of $10^5$ dollars.) Similarly, if our expected cost is small, then we care more about additive regret then about multiplicative regret. (If we pay only $10^{-5}$ dollars, we don't care that we might instead pay $10^{-9}$ dollars.) If we fix $\gamma$ to be a cost we are willing to neglect, then hybrid regret $\text{HR}_\gamma(\sigma,a)$ provides meaningful information regardless of the order of magnitude of the expected cost.

So far, we have defined three flavors of worst-case regret, and posed the problem of minimizing that regret. The solution to our agnostic control problem is given by the following result, proved in \cite{bounded2023}.

\begin{thm}\label{thm: 3}
    Fix $[-a_\emx, a_\emx]$, $q_0$, $T$ (and $\gamma$ if we use hybrid regret). Suppose our \emph{PDE Assumption} is satsified. Then there exist a probability measure $d\eprior(a)$, a finite subset $E \subset [-a_\emx, a_\emx]$, and a strategy $\sigma$, for which the following hold. 
    \begin{enumerate}
        \item[\emph{(I)}] $\sigma$ is the optimal Bayesian strategy for the prior probability distribution $d\eprior$.
        \item[\emph{(II)}] $d\eprior$ is supported in the finite set $E$.
        \item[\emph{(III)}] $E$ is precisely the set of points at which the function $a \mapsto \regret(\sigma, a)$ achieves its maximum on the interval $[-a_\emx, a_\emx]$.
        \item[\emph{(IV)}] $\regret^*(\sigma) \le \regret^*(\sigma')$ for any other strategy $\sigma'$.
    \end{enumerate}
\end{thm}

So, for optimal agnostic control, we should pretend to believe that the unknown $a$ is confined to a finite set $E$ and governed by the probability distribution $d\prior$, even though in fact we know nothing about $a$ except that it lies in $[-a_\mx, a_\mx]$. 

It is easy to see that conditions (I), (II), (III) in Theorem \ref{thm: 3} imply condition (IV) (we give the argument later in this Section). The hard part of Theorem \ref{thm: 3} is the assertion that there exist $d\prior$, $E$, $\sigma$ satisfying (I), (II), (III); we now give an overview of how this is done.

We first prove an analogous result for the setting in which the unknown $a$ is confined to a finite subset $A \subset [-a_\mx, a_\mx]$. Once that's done, we take a sequence of fine nets, e.g., 
\[
A_n = [-a_\mx, a_\mx] \cap 2^{-n} \mathbb{Z}, \; n = 1,2,3,\dots
\]
and deduce Theorem \ref{thm: 3} by applying our result to the $A_n$ and passing to the limit.

We sketch the ideas for finite $A$.

First of all, because we allow strategies to depend on coinflips, it's easy to define intermediate or ``mixed'' strategies between two given strategies $\sigma_0$ and~$\sigma_1$. Given a number $\theta \in [0,1]$, we play strategy $\sigma_1$ with probability $\theta$, and we play instead strategy $\sigma_0$ with probability $1-\theta$. We write $\sigma_\theta$ to denote that mixed strategy. Clearly, we have
\[
\ecost(\sigma_\theta,a) = \theta \ecost(\sigma_1,a) + (1-\theta)\ecost(\sigma_0,a)\;\text{for any}\; a \in \mathbb{R}.
\]

Now let $A \subset [-a_\mx, a_\mx]$ be finite. We associate to any given strategy $\sigma$ its \emph{cost vector}, defined as
\[
\overrightarrow{\ecost}(\sigma) = (\ecost(\sigma, a))_{a \in A} \in \R^A.
\]
Thanks to our discussion of intermediate strategies, the set of all cost vectors of arbitrary strategies is a convex set $\cK \subset \R^A$.

For $\varepsilon > 0$, we call a strategy $\sigma_0$ \emph{$\varepsilon$-efficient} if there is no competing strategy $\sigma'$ such that
\[
\ecost(\sigma', a) < \ecost(\sigma_0,a) - \varepsilon\;\text{for all} \; a \in A.
\]
A simple convexity argument shows that any $\varepsilon$-efficient strategy $\sigma_0$ is within $\varepsilon$ of optimal for some Bayesian prior probability distribution $(p(a))_{a \in A}$ on $A$. To see this, we form the convex set $\cK_-$, consisting of all vectors $(v_a)_{a \in A} \in \R^A$ such that
\[
v_a < \ecost(\sigma_0, a) - \varepsilon\;\text{for all} \; a \in A.
\]

Since $\sigma_0$ is $\varepsilon$-efficient, the convex sets $\cK$ and $\cK_-$ are disjoint, hence there is a nonzero linear functional $\lambda$ on $\R^A$ such that $\lambda(v_-) \le \lambda(v)$ for all $v_- \in \cK_-$, $v \in \cK$. From the functional $\lambda$ we can easily read off a probability distribution $(p(a))_{a \in A}$ on $A$ such that
\[
\sum_{a \in A} p(a) \ecost(\sigma_0, a) \le \sum_{a \in A} p(a) \ecost(\sigma',a) + \varepsilon
\]
for every competing strategy $\sigma'$. 

Thus, as claimed, any $\varepsilon$-efficient strategy is within $\varepsilon$ of best possible for Bayesian control for some prior probability distribution on $A$. Now we are ready for the analogue of Theorem \ref{thm: 3} for finite $A$. The result is as follows. 

\begin{lem}[Agnostic Control Lemma]\label{lem: acl}
    Let $A \subset [-a_\emx, a_\emx]$ be finite, and let $\varepsilon > 0$ be given. Then there exist a subset $A_0 \subset A$, a probability measure $\mu$ on $A_0$, and a strategy $\sigma$ with the following properties.
    \begin{itemize}
        \item $\sigma$ is the optimal Bayesian strategy for the prior $\mu$. 
        \item $\regret(\sigma, a) \le \regret(\sigma, a_0) + \varepsilon$ for all $a \in A$ and $a_0 \in A_0$.
    \end{itemize}
    In particular,
    \[
    |\regret(\sigma, a_0) - \regret(\sigma, a_0')| \le \varepsilon\;\text{for}\; a_0, a_0' \in A.
    \]
\end{lem}

The proof of the Agnostic Control Lemma proceeds by induction on $\#A$, the number of elements of $A$. (So it is essential that the Lemma deals only with finite $A$.)

\underline{In the base case $\#A = 1$}, we have $A = \{a_0\}$ for some $a_0$. We take $A_0 = A$, $\mu = \text{point mass at } a_0$, $\sigma = \text{optimal known-}a \text{ strategy for } a=a_0$. The conclusions of the Lemma are obvious.

\underline{For the induction step}, we fix $k \ge 2$ and suppose our Lemma holds whenever $\# A < k$. We then prove the Lemma for $\#A = k$.

Thus, let $\#A = k$, and let $\varepsilon > 0$. We define suitable small positive numbers
\[
\varepsilon_{4} \ll \varepsilon_3 \ll \dots \ll \varepsilon_0 = \varepsilon.
\]
For $A' \subset [-a_\mx, a_\mx]$ finite, we define
\[
\regret_\mx(\sigma, A') = \max \{ \regret(\sigma,a): a \in A'\}
\]
for any strategy $\sigma$. 

Let $\hat{\sigma}$ be a strategy for which $\regret_\mx(\hat{\sigma}, A)$ is within $\varepsilon_4$ of least possible. Then $\hat{\sigma}$ is $\varepsilon_3$-efficient. Indeed, if any competing strategy $\sigma'$ satisfied
\[
\ecost(\sigma', a) < \ecost(\hat{\sigma},a ) - \varepsilon_3\;\text{for all}\; a \in A,
\]
then $\regret_\mx(\sigma',A)$ would be smaller than $\regret_\mx(\hat{\sigma}, A)$ by more than $\varepsilon_{4}$, contradicting the defining property of $\hat{\sigma}$. Since $\varepsilon_3$-efficient strategies are within $\varepsilon_3$ of best possible for some Bayesian prior, there exists a probability distriubtion $\mu$ on $A$ such that
\begin{equation}\label{eq: l1}
\ecost(\hat{\sigma}, \mu) \le \ecost(\sigma', \mu) + \varepsilon_3    
\end{equation}
for any competing strategy $\sigma'$.

In particular, let $\sigma$ be the optimal Bayesian strategy for the prior $\mu$. Then \eqref{eq: l1} gives
\[
\ecost(\hat{\sigma}, \mu) \le \ecost(\sigma, \mu) + \varepsilon_3.
\]
Theorem \ref{thm: 2}\footnote{Theorem \ref{thm: 2} applies only to tame strategies. In this article, we oversimplify by ignoring that issue. See \cite{bounded2023} for a correct discussion.} therefore implies that
\[
|\ecost(\hat{\sigma}, a) - \ecost(\sigma, a) | \le \varepsilon_3\;\text{for all} \; a \in A,
\]
and therefore
\[
| \regret_\mx(\hat{\sigma}, A) - \regret_\mx(\sigma, A) | \le \varepsilon_2.
\]
Together with the defining property of $\hat{\sigma}$, this shows that
\begin{equation}\label{eq: l2}
    \regret_\mx(\sigma, A) \le \regret_\mx(\sigma', A) + 2 \varepsilon_2.
\end{equation}
for any competing strategy $\sigma'$.

It may happen that
\begin{equation}\label{eq: l3}
    \regret(\sigma, a) \ge \regret_\mx(\sigma, A) - \varepsilon_1\;\text{for all}\; a \in A.
\end{equation}
In that case, we have
\[
\regret_\mx(\sigma,A) - \varepsilon_4 \le \regret(\sigma, a) \le \regret_\mx(\sigma, A)\;\text{for all}\; a \in A,
\]
so the conclusions of our lemma hold for the above $\mu$, $\sigma$ with $A_0 = A$. Hence, we may assume that \eqref{eq: l3} is false.

We set
\[
A_0 = \{a \in A: \regret(\sigma, A) \ge \regret_\mx(\sigma, a) - \varepsilon_1\}.
\]
Since \eqref{eq: l3} is false, we have $\# A_0 < \#A = k$, hence, by our induction hypothesis, Lemma \ref{lem: acl} applies to $A_0$.

Thus, there exist a subset $A_{00} \subset A_0$, a probability measure $\mu_0$ on $A_{00}$, and a strategy $\sigma_0$, such that
\begin{itemize}
    \item $\sigma_0$ is the optimal Bayesian strategy for the prior $\mu_0$, and
    \item $\regret(\sigma_0, a) \le \regret(\sigma_0, a_0) + \varepsilon_4$ for all $a \in A_0$, $a_0 \in A_{00}$.
\end{itemize}
We then show that the conclusions of Lemma \ref{lem: acl} hold, with $A_{00}$, $\mu_0$, $\sigma_0$ in place of $A_0$, $\mu$, $\sigma$. This completes our induction on $\# A$, proving Lemma \ref{lem: acl}.

Once we have established Lemma \ref{lem: acl}, we can easily pass from the finite sets $A_n = [-a_\mx, a_\mx] \cap 2^{-n} \mathbb{Z}$ to the full interval $[-a_\mx, a_\mx]$ by a weak compactness argument. This proves conclusions (I), (II), (III) of Theorem \ref{thm: 3} except for the finiteness of the set $E$ on which the function
\[
[-a_\mx, a_\mx] \ni a \mapsto\regret(\sigma,a)
\]
takes its maximum.

To see that $E$ is finite, we examine the function
\[
F: \R \ni a \mapsto \regret(\sigma,a ).
\]
We prove that $F$ is real-analytic and grows exponentially fast as $a \rightarrow + \infty$. Consequently, $F|_{[-a_\mx, a_\mx]}$ is a nonconstant real-analytic function, which can therefore achieve its maximum at only finitely many points. Thus, (I), (II), and (III) hold with $E$ finite.

It remains only to deduce conclusion (IV) from (I), (II), (III). Let $d\prior$, $\sigma$, $E$ be as in (I), (II), (III) of Theorem \ref{thm: 3}. Since $\sigma$ is the optimal Bayesian strategy for $d\prior$ (by (I)), and since $d\prior$ is supported on the finite set $E$ (by (II)), we have for any other strategy $\sigma'$ that
\[
\ecost(\sigma, a_0) \le \ecost(\sigma', a_0) \;\text{for some} \; a_0 \in E.
\]
In particular, we have
\[
\regret(\sigma, a_0) \le \regret(\sigma' , a_0)\;\text{for some}\; a_0 \in E.
\]
Combining this with (III), we see that for any $a \in [-a_\mx, a_\mx]$ we have
\[
\regret(\sigma, a) \le \regret(\sigma',a_0).
\]
Therefore (I), (II), (III) of Theorem \ref{thm: 3} easily imply (IV).

This concludes our discussion of agnostic control for bounded $a$; for details, see \cite{bounded2023}. Finally, we pass to the most general case.

\subsection*{Variant IV: Fully Agnostic Control}

Finally, we make no assumption whatever regarding the unknown $a$. Our $a$ may be any real number, and we are not given a Bayesian prior distribution for it. If $a$ is large positive, then the system is highly unstable. Our goal is again to minimize worst-case regret, defined as in the previous section, except now the supremum is taken over all $a \in \R$. We confine ourselves to hybrid regret.

We now denote the hybrid regret of a strategy $\sigma$ by $\text{HR}_\gamma(\sigma, a; q_0, T)$, to make explicit the r\^{o}le of the starting position $q_0$ and time horizon $T$. Thus, for fixed $\gamma$, $q_0$, $T$, we are trying to minimize
\[
\text{HR}_\gamma^*(\sigma; q_0, T) = \sup_{a \in \R} \text{HR}_\gamma(\sigma; a, q_0, T).
\]
We remark that this sup may be infinite.

We strengthen our \emph{PDE Assumption} by assuming also that the constant $C_\tame$ in \eqref{eq: 17} grows at most as a power of $a_\mx$ when $a_\mx \gg 1$, i.e., we assume that \eqref{eq: 11}--\eqref{eq: 16} hold and that there exists an integer $n_0$ for which
\begin{equation}\label{eq: 18}
|\tilde{u}| \le C_0\cdot [1+a_\mx^{n_0}] \cdot [1 + |q|] \;\text{for all}\; (q, t, \zeta_1, \zeta_2)
\end{equation}
(recall that $a_\mx>0$). This seems plausible; we have argued that most likely $C_\tame = O(a_\mx)$ (see \eqref{eq: 4.5}).

The main result of our paper \cite{almostoptimal2023} is that, with negligible increase in regret, we can reduce matters to agnostic control for bounded $a$. Specifically, we prove the following Theorem.

\begin{thm}\label{thm: 4}
    Fix a time horizon $T$, a nonzero starting position $q_0$, and constants $C_0, n_0$ (to be used in the estimate \eqref{eq: 18}). Then given $\varepsilon>0$ there exists $a_\emx>0$ for which the following holds.
    
    Let $\sigma$ be a strategy for the starting position $q_0$ and time horizon $T+\varepsilon$. Suppose $\sigma$ satisfies estimate \eqref{eq: 18} for $a_\emx$ and the given $C_0, n_0$. 
    
    Then there exists a strategy $\sigma_*$ for the starting position $q_0$ and time horizon $T$, satisfying the following estimates.
    \begin{itemize}
        \item[\emph{(A)}] For $a \in [-a_\emx,a_\emx]$ we have
        \begin{multline*}
        \ecost(\sigma_*,a;T,q_0)\\ \le \varepsilon + (1+\varepsilon)\cdot \sup\big\{\ecost(\sigma,a';T+\varepsilon,q_0) : |a'- a| \le \varepsilon |a|\big\}.
        \end{multline*}
        \item[\emph{(B)}] For $a \notin [-a_\emx,a_\emx]$ we have
        \[
        \ecost(\sigma_*,a;T,q_0) \le \varepsilon + (1+\varepsilon)\cdot \ecost(\sigma_{\emph{opt}}(a),a;T,q_0).
        \]
    \end{itemize}
\end{thm}

So, if $a \in [-a_\mx,a_\mx]$, then $\sigma_*$ performs almost as well as $\sigma$; and if $a \notin [-a_\mx,a_\mx]$, then $\sigma_*$ performs almost as well as the optimal known-$a$ strategy $\sigma_\opt(a)$.

Using Theorem \ref{thm: 4}, we construct strategies $\sigma$ that come arbitrarily close to minimizing worst-case hybrid regret. Assume that we are given constants $\gamma, T, q_0, C_0, m_0$ as in Theorem \ref{thm: 4}. We let $\varepsilon>0$ be given and we take $a_\mx$ to be a large enough positive real number (depending on $\varepsilon$ as well as the constants above).

We let $\sigma_0$ be the optimal agnostic control strategy for worst-case hybrid regret with starting position $q_0$ and time horizon $T+\varepsilon$, and with $a$ confined to the interval $[-(1+\varepsilon)a_\mx, (1+\varepsilon)a_\mx]$. (Of course, this is Variant III above). We assume \eqref{eq: 18} holds for $\sigma_0$.

Applying Theorem \ref{thm: 4} to $\sigma_0$, we obtain a strategy $\sigma_{\text{Ag}}$ for time horizon $T$ so that: 
\begin{itemize}
    \item For $a \in [-a_\mx, a_\mx]$, the strategy $\sigma_{\text{Ag}}$ performs only slightly worse than the worst-case performance of the strategy $\sigma_0$ on the slightly larger interval $[-(1+\varepsilon)a_\mx, (1+\varepsilon)a_\mx]$.
    \item For $a \notin [-a_\mx, a_\mx]$, the strategy $\sigma_{\text{Ag}}$ performs only slightly worse than the optimal known-$a$ strategy $\sigma_\opt(a)$.
\end{itemize}
From this, it's easy to deduce that the worst-case hybrid regret of the strategy $\sigma_{\text{Ag}}$ (for fully agnostic control, i.e., with $a \in \R$) is at most $O(\varepsilon)$ percent worse than that of $\sigma_0$ (for agnostic control with $a$ confined to $[-(1+\varepsilon)a_\mx, (1+\varepsilon)a_\mx]$). The worst-case hybrid regret of the optimal strategy $\sigma_0$ on the interval $[-(1+\varepsilon)a_\mx, (1+\varepsilon)a_\mx]$ is, of course, bounded above by the worst-case hybrid regret of \emph{any} strategy $\sigma$ for fully agnostic control (i.e., with $a \in \R$). Consequently, we have
\[
\text{HR}_\gamma^*(\sigma_\ag; q_0, T) \le (1+ C\varepsilon) \cdot \text{HR}_\gamma^*(\sigma; q_0 , T+\varepsilon)
\]
for any competing strategy $\sigma$.

Thus, building on our solution for the control problem in Variant III, we have produced an almost optimal strategy for fully agnostic control. For a more detailed overview of the proof of Theorem \ref{thm: 4}, we refer the reader to the introduction of \cite{almostoptimal2023}.

\section*{A Future Direction}

In \cite{bounded2023}, we discuss several unsolved problems suggested by our work in \cite{almostoptimal2023,bounded2023}. Here, we discuss one of those unsolved problems in more detail. Specifically, we speculate briefly on a particular model problem in which we don't know \emph{a priori} what our control does.

Consider a particle governed by the stochastic ODE
\begin{equation}\label{eq: fd 1}
    dq(t) = au(t)dt + dW(t), \qquad q(0) = 0.
\end{equation}
As usual, $q(t)$ denotes position, $u(t)$ is our control, $W(t)$ is Brownian motion, and we incur a cost
\[
\int_0^T \{ (q(t))^2 + (u(t))^2\} \ dt.
\]
We would like to understand optimal agnostic control for this system, i.e., we'd like to find strategies that minimize worst-case regret. In analogy with our work on the system \eqref{eq: intro 1}, we first attempt to understand optimal Bayesian control.

In the simplest case of Bayesian control, suppose we know \emph{a priori} that $a=1$ or $a=-1$, each with probability $1/2$.

We write $\ecost(\sigma)$ to denote the expected cost incurred by executing a strategy $\sigma$, and we set
\begin{equation}\label{eq: fd 2}
    \ecost^* = \inf \{ \ecost(\sigma): \;\text{All strategies } \sigma\}.
\end{equation}
For this simple problem, we make the following conjectures.
\begin{itemize}
    \item The infimum in \eqref{eq: fd 2} is not achieved by any strategy $\sigma$, because there is a regime in which we would like to set $u(t) = \pm \infty$, in order to gain instant information about $a$.
    \item A nearly optimal strategy will determine $u(t)$ as a function of position $q(t)$, time $t$, and $p(t) = \text{posterior probability that} \;a = + 1$, given history up to time $t$. Thus, $u(t) = \tilde{u}(q(t), t, p(t))$ for a function $\tilde{u}(q,t,p)$ on the ``state space'' $\Omega = \R \times [0,T] \times [0,1]$.
    \item The state space $\Omega$ is partitioned into two regimes $\Omega_0$ and $\Omega_1$. In $\Omega_0$, we would like to set $\tilde{u} = \pm \infty$, so we set $\tilde{u} = \mathcal{U}$, a large positive number.\footnote{We could just as well set $\tilde{u} = - \mathcal{U}$.} In $\Omega_1$, we take $\tilde{u}$ to be a solution of a relevant Bellman equation. A free boundary condition determines how we partition $\Omega$ into $\Omega_0$ and $\Omega_1$.
    \item As $\mathcal{U} \rightarrow \infty$, such strategies approach optimality. Perhaps one should define strategies in a way that allows $u = \pm \infty$. If so, this had better be done carefully.
\end{itemize}

We emphasize that the above are speculations---we have no rigorous results on optimal agnostic control for the system \eqref{eq: fd 1}. We remark, however, that in \cite{carruth2023bounded} the first-named author has found a strategy that achieves bounded multiplicative regret for a more general system than \eqref{eq: fd 1}.\footnote{In fact, the results of \cite{carruth2023bounded} assume that the starting position $q_0$ satisfies $|q_0| \ge 1$. After an easy modification, however, the strategy defined in that paper for the system \eqref{eq: fd 1} achieves bounded regret for arbitrary $q_0 \in \R$.}

Clearly, there is much to be done before we can claim to understand agnostic control theory.

\bibliographystyle{plain}
\bibliography{ref}
\end{document}